\magnification = \magstep1



\font\Bbb=msbm10
\font\sll  = cmss10

\font\bfi  = cmmib9

\font\msbma = msbm9

\def\BbbN{\hbox{\Bbb N}}
\def\BbbQ{\hbox{\Bbb Q}}
\def\BbbR{\hbox{\Bbb R}}
\def\BbbZ{\hbox{\Bbb Z}}
\def\setmin{\mathop{\raise.1ex\hbox{\msbma\char31}}}

\baselineskip=18pt

\hsize13.2truecm 

\vsize21.8truecm \vbox to 2.1truecm{}

\bf
\centerline{\bf ON THE MINIMAL IDEMPOTENTS OF}
\centerline{\bf TWISTED GROUP ALGEBRAS OF CYCLIC 2-GROUPS}

\vskip2pc\vskip9pt\bf
\centerline{{\bfi Jordan} {\bfi J}.\thinspace{\bfi Epitropov},
	    {\bfi Todor} {\bfi Zh}.\thinspace{\bfi Mollov},
	    {\bfi Nako}   {\bfi A}.\thinspace{\bfi Nachev} \footnote *
       {\baselineskip=11pt\sll
	   Supported by the National Scientific Fund of the Ministry
	   of Education, Science  and Technologies of Bulgaria under
	   Contract MM431/94.
       }}
\vskip1pc\vskip3pt
\centerline{{\bfi Department} {\bfi of} {\bfi Algebra},
	    {\bfi University} {\bfi of} {\bfi P\/lovdiv},}
\centerline{4000 {\bfi P\/lovdiv}, {\bfi Bulgaria}}

\vskip4pc\vskip20pt
\leftskip.96truecm\hsize12.24truecm 
\sll
\baselineskip=11pt

{\bf Abstract.} For a field $K$ of the second kind with respect
to 2 and of characteristic different from 2, we consider the
decomposition of the binomials $x^{2^n} - a$ into a product of
irreducible factors over $K$ and find the explicit
form of the minimal idempotents of the twisted group algebra
$K^t\langle g\rangle$ of a cyclic 2-group $\langle g\rangle$ over
$K$.

\vskip2pc\leftskip0truecm
\hsize13.2truecm
\rm\baselineskip=18pt

\centerline{\bf Introduction}
\medskip

The starting point of a lot of investigations on twisted group
algebras is to find the minimal idempotents of the twisted group
algebra $K^tG$ where $G$ is a cyclic group and $K$ is a field. When
$\langle g\rangle$ is a cyclic $p$-group, $p$ is an odd prime and
$K$ is a field of characteristic different from $p$, Nachev and
Mollov [3] have found the explicit form of the minimal idempotents
of $K^t\langle g\rangle$. For $p = 2$
  additional difficulties arise which are connected
 with the decomposition of the
polynomial $x^{2^n} - a$ into irreducible factors over the field
$K$. The purpose of the present paper is to find the explicit form
of the minimal idempotents of $K^t\langle g\rangle$ when
$\langle g\rangle$ is a cyclic 2-group and $K$ is a field of the
second kind with respect to 2 (and of characteristic different from
2). We shall mention that the semisiplicity of the twisted group
algebra $K^t\langle g\rangle$ in this case is a well known fact
(see e.g. [1] or [4]).

The paper is organized as follows. In Section 1 we give some
notation, definitions and preliminary results. Section 2 deals with
the decomposition of an arbitrary polynomial $x^{2^n} - a$ in
irreducible factors over a field $K$ of the second kind with
respect to 2. In Section 3 we find the explicit form of the minimal
idempotents of the twisted group algebra $K^t\langle g\rangle$ of
an arbitrary cyclic 2-group $\langle g\rangle$ over a field $K$ of
the second kind with respect to 2.
\vskip 1.5truecm

\centerline{\bf 1. Notation, Definitions and Preliminary Results}
\bigskip

Let $p$ be a prime, let $K$ be a field of characteristic different
from $p$ and let $\bar K$ be the algebraic closure of $K$. We
denote by $\varepsilon_n$ a $p^n$-th primitive root of 1 in
$\bar K$. The field $K$ is called a field of the first kind with
respect to $p$ [5, p. 684], if
$K(\varepsilon_i)\not= K(\varepsilon_2)$ for some $i > 2$.
Otherwise $K$ is called a field of the second kind with respect to
$p$. An equivalent definition is the following. If the Sylow
$p$-subgroup $K(\varepsilon_2)_p$ of the multiplicative group
$K(\varepsilon_2)^{\ast}$ is cyclic, then $K$ is of the first kind
with respect to $p$ and if $K(\varepsilon_2)_p$ is the quasicyclic
group $\BbbZ(p^{\infty})$ then $K$ is a field of the second
kind with resect to $p$. Typical examples of fields of the
first and of the second kind with respect to any prime $p$
are $\BbbQ$ and $\BbbR$, respectively. Let
$K^{\ast}$ be the multiplicative group of $K$ and let $G$ be a
multipliucative group. A twisted group algebra $K^tG$ of $G$ over
$K$  [5, p. 13] is an associative $K$-algebra with basis
$\{\bar x\mid x\in G\}$ and with multiplication defined on the basis by
$$
\bar x\bar y = \gamma(x,y)\overline{xy},\, \gamma(x,y)\in K^{\ast}.
$$
We denote by $KG$ the ordinary group algebra of $G$ over $K$. It is
well known that if $G=\langle g\rangle$ is a cyclic group then
$K^t\langle g\rangle$ is a commutatitve algebra. If the group
$\langle g\rangle$ is of order $n$ then $\bar g^n = a$ for some
$a\in K^{\ast}$. Obviously this equality determines the twisted
group algebra $K^t\langle g\rangle$. We define
$$
K^n =\{a^n\mid a\in K\},\,n\in\BbbN.
$$
Clearly $K^n$ is closed with respect to the multiplication and
especially $K^{p^n}\supseteq K^{p^{n+1}}$.

The following definition is well known [6, \S 44, p. 142].
\medskip

{\bf Definition 1.1.} Let $\alpha$ belong to the finite extension
$L$ of the field $K$ and let $f(x) = x^n + a_1x^{n-1}+\ldots+a_n$
be the minimal polynomial of $\alpha$ over $K$. The element
$$
N(\alpha)=(-1)^na_n^{(L:K)/n}
$$
is called the norm of $\alpha$ (in $L$ over $K$).
\medskip

Clearly $N(\alpha)= (\alpha_1\ldots\alpha_n)^{(L:K)/n}$, where
$\alpha_1,\ldots,\alpha_n$ are the zeros of the polynomial $f(x)$,
i.e. all the conjugate elements of $\alpha$ over $K$. It is well
known that $N(\alpha\beta)=N(\alpha)N(\beta)$ for
$\alpha,\beta\in L$ and if $a\in K$ then $N(a) = a^{(L:K)}$.

In order to see that some special binomials are indecomposable over
the field $K$ we shall use the following theorem.
\medskip

{\bf Theorem 1.2 {\rm [2, Theorem 16.6, p. 225]}.} {\sl The
binomial $x^n - a$, $a\in K$, is irreducible over $K$ if and only
if $a\not\in K^p$ for all primes $p$ dividing $n$ and
$a\not\in -4K^4$ whenever $4\vert n$.}
\medskip

Allover in this paper the base field $K$ will be of the second kind
with respect to 2 and of characteristic different from 2. If
$a\in K$ and $n$ is a fixed positive integer, then we denote by
$H_n(a)$ the greatest integer $s$ in the interval $[0,n]$ such that
$a\in K(\varepsilon_2)^{2^s}$. As usually we assume that 0 and 1
are the trivial idempotents of an algebra.
\vskip 1.truecm

\centerline{\bf 2. Decomposition of Special Binomials over a Field}
\centerline{\bf in Irreducible Factors}
\medskip

If it is not explicitly stated, in this section we assume that
$K\not=K(\varepsilon_2)$, i.e. $-1\not\in K^2$.
\medskip

{\bf Lemma 2.1.} {\sl For every $n\geq 2$ the only conjugated over
$K$ element $\bar\varepsilon_n$ of $\varepsilon_n$ is
$\varepsilon_n^{-1}$.}

{\sl Proof.} Since $\varepsilon_n$ is a root of the equation
$x^{2^{n-1}} + 1 = 0$ over $K$, we obtain that $\bar\varepsilon_n$
is a root of the same equation. Hence
$\bar\varepsilon_n\in \langle\varepsilon_n\rangle$ and
$\varepsilon_n\bar\varepsilon_n\in
K\cap\langle\varepsilon_n\rangle=\{-1,1\}$,
i.e. $\bar\varepsilon_n=t\varepsilon_n^{-1}$, where $t\in\{-1,1\}$.
Besides
$$
\varepsilon_n + \varepsilon_n^{-1} =
(\varepsilon_{n+1} + t\varepsilon_{n+1}^{-1})^2 - 2t =
(\varepsilon_{n+1} + \bar\varepsilon_{n+1})^2 - 2t \in K.
$$
Therefore $\varepsilon_n + \varepsilon_n^{-1}\in K$ and
$\varepsilon_n\varepsilon_n^{-1} = 1\in K$, which shows that
$\bar\varepsilon_n = \varepsilon_n^{-1}$.
\medskip

{\bf Lemma 2.2.} {\sl Let $\alpha\in K(\varepsilon_2)$. Then the
following conditions are equivalent.}

(i) $\alpha\in K$;

(ii) {\sl $\bar \alpha=\alpha$, where $\bar \alpha$ is the conjugated of $\alpha$;}

(iii) $N(\alpha) = \alpha^2$.

{\sl Proof.} The equivalence of (i) and (ii) is obvious. The
equivalence of (ii) and (iii) follows from the equality
$\alpha^2 = N(\alpha) = \alpha\bar \alpha$.
\medskip

{\bf Lemma 2.3.} {\sl The polynomial $f(x) = x^{2^n} - a$,
$a\in K$, $n\geq 2$, is irreducible over $K$ if and only if
$a\not\in K^2\cup(-K^4)$.}

  Really since $4 = (\varepsilon_3 + \varepsilon_3^{-1})^4 =
(\varepsilon_3 + \bar\varepsilon_3)^4\in K^4$, then
  $4K^4 = K^4$. Then the proof of the lemma follows from Theorem 1.2.

\medskip

{\bf Lemma 2.4.} {\sl Let $\alpha\in K(\varepsilon_2)$. Then
$\alpha^{2^n}\in K$ if and only if $\alpha=a\varepsilon_{n+1}^l$,
$a\in K$, $l\in\BbbZ$.}

{\sl Proof.} Let $\alpha^{2^n}\in K$. By Lemma 2.2 we see that
$(N(\alpha))^{2^n} = N(\alpha^{2^n}) = \alpha^{2^{n+1}}$. Therefore
we obtain $N(\alpha) = \alpha^2\varepsilon_n^{-l} =
(\alpha\varepsilon_{n+1}^{-l})^2$ for some $l\in\BbbZ$. On
the other hand, by Lemma 2.1 we have that
$N(\varepsilon_{n+1}^{-l}) = 1$. Hence
$N(\alpha\varepsilon_{n+1}^{-l}) = N(\alpha) =
(\alpha\varepsilon_{n+1}^{-l})^2$ and by Lemma 2.2 it follows that
$\alpha\varepsilon_{n+1}^{-l}\in K$, i.e. $\alpha =
a\varepsilon_{n+1}^l$. The inverse statement of the lemma is
obvious.
\medskip

{\bf Lemma 2.5.} {\sl For every $n\in\BbbN$ it holds}
$$
K\cap K(\varepsilon_2)^{2^n} = K^{2^n}\cup (-K^{2^n}).
$$

{\sl Proof.} Let $u\in K\cap K(\varepsilon_2)^{2^n}$. Then $u=\alpha^{2^n},\
\alpha\in K(\varepsilon_2)$ and, by Lemma 2.4, $\alpha=b\varepsilon_{n+1}^l,\
b\in K, l\in \BbbZ$. Therefore
$$
  u = \alpha^{2^n} = b^{2^n}(-1)^l \in K^{2^n} \cup (-K^{2^n}),
$$
i.e.
$$
K\cap K(\varepsilon_2)^{2^n} \subseteq K^{2^n} \cup (-K^{2^n}).
$$
The opposite inclusion is obvious.
\medskip

{\bf Lemma 2.6.} {\sl The equality $K^{2^n}\cap (-K^{2^m}) = 0$
holds for arbitrary $m,n\in\BbbN$.}

{\sl Proof}. The equality $K^2 \cap (-K^2) = 0$ follows from the fact that
$\varepsilon_2\not\in K$.
Moreover $K^{2^n}\cap (-K^{2^m})\subseteq K^2\cap (-K^2) = 0$.
\medskip

{\bf Lemma 2.7.} {\sl The polynomial $f(x) = x^{2^n} - \alpha$,
$\alpha\in K(\varepsilon_2)$, $n\in\BbbN$, is irreducible
over $K(\varepsilon_2)$ if and only if $\alpha\not\in
K(\varepsilon_2)^2$.}

For $n\ge 2$ the proof follows from Lemma 2.3 applied to $K(\varepsilon_2)$
bearing in mind that
$-K(\varepsilon_2)^4\subseteq K(\varepsilon_2)^2$.
\medskip

 For $n\ge1$ the proof follows from Theorem 1.2 applied to the field
$K(\varepsilon_2)$.

Let $a\in K^{\ast}$ and $n\in\BbbN$. We denote by
$H_n(a)$ the greatest integer $s$ in the interval $[0,n]$ such that
$a\in K(\varepsilon_2)^{2^s}$. Lemma 2.5 gives immediately that the
integer $H_n(a)$ coincides with the greatest integer $s\in [0,n]$
such that $a\in K^{2^s}\cup (-K^{2^s})$. If $a\in K^{2^s}$, then we
call the element $a$ of {\sl the first kind} and if $a\in -K^{2^s}$,
then we call $a$ an element of {\sl the second kind}. Clearly the
integer $H_n(a)$ always exists and is uniquely determined, by Lemma
2.6 the kind of $a$ is also completely determined. We call the
integer $H_n(a)$ the $n$-{\sl height} of $a$ in $K^{\ast}$.
\medskip

{\bf Theorem 2.8.} {\sl Let $K$ be a field of the second kind with
respect to $2$ and of characteristic different from $2$, let
$f(x)=x^{2^n} - a$ be a polynomial over $K$, $a\not= 0$,
$n\in\BbbN$, and let $H_n(a) = s$. Then}

1) {\sl If $K = K(\varepsilon_2)$, then}
$$
f(x) = \prod_{i=0}^{2^s-1}(x^{2^{n-s}} -b\varepsilon_s^i),\,
a = b^{2^s},\,b\in K.
\leqno (1)
$$
and the factors of (1) are irreducible polynomials over $K$.

2) {\sl if $K\not= K(\varepsilon_2)$  and} 2.1) $s = 0$ {\sl or}
2.2) $s = 1$ {\sl and the element $a$ is of the second kind, then
$f(x)$ is irreducible over $K$.

In the other cases $f(x)$ is
decomposed in irreducible factors over $K$ in the following way.}

3) {\sl If $K\not= K(\varepsilon_2)$, $s\geq 1$, and $a$ is of the first
kind, then}
$$
f(x) = (x^{2^{n-s}} - b)(x^{2^{n-s}} + b)
\prod_{i=1}^{2^{s-1}-1}
[x^{2^{n-s+1}} - (\varepsilon_s^i + \varepsilon_s^{-i})bx^{2^{n-s}}
+ b^2],
\leqno (2)
$$
$a = b^{2^s},\,b\in K.$

4) {\sl If $K\not= K(\varepsilon_2)$, $s\geq 2$, and $a$ is of
the second kind, then}
$$
f(x) = \prod_{i=0}^{2^{s-1}-1}
[x^{2^{n-s+1}} - (\varepsilon_s^i + \varepsilon_s^{-i-1})
\varepsilon_{s+1}bx^{2^{n-s}} + b^2],\,
a = -b^{2^s}, \leqno (3)
$$
$b\in K$.

{\sl Proof.} 1) Let $K = K(\varepsilon_2)$. Obviously (1) is a
decomposition of $f(x)$ over $K$. If $s < n$, then
$b\varepsilon_s^i\not\in K^2,\ -K^4\subseteq K^2$
and by Lemma 2.7 the factors of (1)
are irreducible polynomials over $K$. The case $s = n$ is
trivial.

2) Let $K \not= K(\varepsilon_2)$ and 2.1) $s = 0$ or 2.2) $s = 1$
and $a$ is of the second kind.

  If $n = 1$, then in the case 2.1) we
have $a\not\in K^2$ and in the case 2.2) by Lemma 2.6 we have again
$a\not\in K^2$. Now Theorem 1.2 gives the irreducibility of $f(x)$
over $K$.

  Let $n\ge 2$. If $s=0$, then $a\not\in K(\varepsilon_2)^2
\supset (-K^4)$, i.e. $a\not\in K^2 \cup (-K^4)$. Therefore, by Lemma 2.3,
the polynomial $f(x)$ is irreducible over $K$. Let $s=1$ and
the element $a$ is of the second kind. Then, by the definition of $a$,
it follows that $a\in (-K^2)\setmin (-K^4)$ and, by Lemma 2.6, $a\not\in
K^2 \cup (-K^4)$. Thus, again by Lemma 2.3, $f(x)$ is irreducible over $K$.

3) Let $K \not= K(\varepsilon_2)$, $s\geq 1$, and let $a$ be of
the first kind. Then $a=b^{2^s}$, $b\in K$. We obtain the decomposition
(1) of the polynomial $f(x)$ over $K(\varepsilon_2)$. We shall show
that all the factors of $f(x)$ in (1) are irreducible over
$K(\varepsilon_2)$. Really, if $s < n$, then by the definition of
$H_n(a)$ we obtain that $b\not\in K(\varepsilon_2)^2$ and, since
$\varepsilon_s^i = \varepsilon_{s+1}^{2i} \in K(\varepsilon_2)^2$, it follows that
$b\varepsilon_s^i\not\in K(\varepsilon_2)^2$. Now, by Lemma 2.7 the
factors in (1) are irreducible over $K(\varepsilon_2)$. For
$s=n$ the factors of (1) are of the first degree and also are
irreducible over $K(\varepsilon_2)$. Grouping and multiplying the
conjugated over $K$ factors in (1), by Lemma 2.1, we obtain that
the factors of
(2) are with coefficients from $K$. Their irreducibility over $K$ follows
from the irreducibility of the factors of (1) over $K(\varepsilon_2)$.

4) Let $K\not= K(\varepsilon_2)$, $s\geq 2$, and let $a$ be of
the second kind. Then $a = -b^{2^s}$, $b\in K$, and we obtain the
following decomposition of the polynomial $f(x)$ over
$K(\varepsilon_2)$
$$
f(x) = \prod_{i=0}^{2^s-1}(x^{2^{n-s}} -
b\varepsilon_{s+1}^{\phantom{i}}\varepsilon_s^i).
\leqno (4)
$$
As in the case 3) we see that (4) is a decomposition of $f(x)$ in
irreducible factors over $K(\varepsilon_2)$ and (4), in view of Lemma 2.1,
gives the
decomposition (3) of $f(x)$ in irreducible factors over $K$.
\vskip 1.5truecm

\centerline{\bf 3. Minimal Idempotents of Twisted Group Algebras}
\centerline{\bf of Cyclic 2-Groups}
\bigskip

{\bf Theorem 3.1.} {\sl Let $K$ be a field of the second kind with
respect to $2$ and of characteristic different from 2 and let
$\langle g\rangle$ be a cyclic group of
order $2^n$. Let the twisted group algebra $K^t\langle g\rangle$ be
defined by the equality $\bar g^{2^n} = a$, $a\in K^{\ast}$, and let
$H_n(a) = s$. Then the minimal idempotents $e_i$ and $f_k$ of
$K^t\langle g\rangle$ are the following.}

1) {\sl If $K = K(\varepsilon_2)$ then
$$
e_i = {1\over 2^s}\sum_{j=0}^{2^s-1}
\varepsilon_s^{-ij}b^{-j}\bar g^{2^{n-s}j},\,
i = 0,1,\ldots,2^s-1,
\leqno (1)
$$
where $b^{2^s} = a$, $b\in K$.}

2) {\sl If $K\not= K(\varepsilon_2)$ and 2.1) $s = 0$ or 2.2) $s = 1$ and
the element $a$ is of the second kind then the only minimal idempotent of
$K^t\langle g\rangle$ is the unity.}

3) {\sl If $K\not= K(\varepsilon_2)$, $s \geq 1$, and $a$ is of
the first kind then
$$
\eqalign{
f_k &= {1\over 2^s}\sum_{j=0}^{2^s-1}
\delta_{kj}b^{-j}\bar g^{2^{n-s}j},\, k = 1,2;\,\delta_{1j}=1,\,
\delta_{2j} = (-1)^j,\cr
e_i &= {1\over 2^s}\sum_{j=0}^{2^s-1}
(\varepsilon_s^{ij} + \varepsilon_s^{-ij})
b^{-j}\bar g^{2^{n-s}j},\,
i = 1,2,\ldots,2^{s-1}-1,
} \leqno (2)
$$
where $b^{2^s} = a$, $b\in K$.}

4) {\sl If $K\not= K(\varepsilon_2)$, $s\geq 2$, and $a$ is of the second kind
then
$$
e_i = {1\over 2^s}\sum_{j=0}^{2^s-1}
(\varepsilon_s^{ij+j} +\varepsilon_s^{-ij})\varepsilon_{s+1}^{-j}
b^{-j}\bar g^{2^{n-s}j},\,
i=0,1,\ldots,2^{s-1} - 1,
\leqno (3)
$$
where $-b^{2^s} = a$, $b\in K$.}

{\sl Proof.} Let $L$ be the splitting field of the polynomial
$f(x) = x^{2^n} - a$ over $K$. It is well known, that the minimal idempotents
of the group algebra $L\langle g\rangle$ are
$$
   e_{\beta} = {1\over 2^n}\sum_{i=0}^{2^n-1}
\beta^{-i} g^i,
\leqno (4)
$$
where $\beta$ runs on the zeros of the polynomial $x^{2^n}-1$. Since $L$
is a splitting field of the polynomial $f(x)$, then there exists $\gamma\in L$,
such that $\gamma^{2^n}=a$. Then the equality $\bar g^{2^n}=a$ implies
$(\gamma^{-1}\bar g)^{2^n} =1$ and therefore the cyclic group $\langle\gamma^{-1}
\bar g\rangle$ is a group basis of the twisted group algebra $L^t\langle g\rangle$,
i.e. $L^t\langle g\rangle$ coincides with the group algebra $L\langle\gamma^{-1}
\bar g\rangle$. Hence the minimal idempotents of $L^t\langle g\rangle$ are
obtained from (4) by replacing $g$ with $\gamma^{-1}\bar g$. So we obtain
that the minimal idempotents of $L^t\langle g\rangle$ will be of the form
$$
   e_{\beta} = {1\over 2^n}\sum_{i=0}^{2^n-1}
(\beta\gamma)^{-i} \bar g^{i},
$$
When $\beta$ runs on the zeros of $x^{2^n}-1$, $\beta\gamma$ will run on
the zeros of $f(x)$. Therefore we can set $\beta\gamma=\alpha$ and the
minimal idempotents of the twisted group algebra $L^t\langle g\rangle$
will be of the form
$$
e_{\alpha} = {1\over 2^n}\sum_{j=0}^{2^n-1}
\alpha^{-j}\bar g^j,
$$
where $\alpha$ runs on the zeros of $f(x)$. From here we shall obtain the
minimal idempotents of $K(\varepsilon_2)^t\langle g\rangle$ summing
the minimal idempotents $e_{\alpha}$, where $\alpha$ runs on
the zeros of
an arbitrary fixed irreducible factor of $f(x)$ over $K(\varepsilon_2)$.
By Theorem 2.8 the irreducible factors of $f(x)$ over $K(\varepsilon_2)$
are of the form   
$\varphi_i(x) = x^{2^{n-s}} - \lambda\varepsilon_s^i$,
$i = 0,1,\ldots,2^s - 1$, where $\lambda^{2^s} = a$, $\lambda\in K(\varepsilon_2)$. Let
$\mu^{2^{n-s}} = \lambda$, $\mu\in L$. The zeros of $\varphi_i(x)$ are
$\mu\varepsilon_n^i\varepsilon_{n-s}^r$, $r = 0,1,\ldots,2^{n-s}-1$,
and the idempotents $e_i$ of $K(\varepsilon_2)^t\langle g\rangle$ are
$$
e_i = {1\over 2^n}\sum_{j=0}^{2^n-1}
\left(\sum_{r=0}^{2^{n-s}-1}
\varepsilon_{n-s}^{-rj}\right)
\mu^{-j}\varepsilon_n^{-ij}\bar g^j.
$$
For a fixed $j$ not divisible by $2^{n-s}$ the sum in the brackets is equal
to zero. Therefore if we replace $j$ by $j2^{n-s}$,
$j=0,1,\ldots,2^{n-s}-1$, we shall obtain
$$
e_i = {1\over 2^s}\sum_{j=0}^{2^s-1}
\varepsilon_s^{-ij}\lambda^{-j}\bar g^{j2^{n-s}},\,
i = 0,1,\ldots,2^s - 1.
\leqno (5)
$$
Now in the case 1) for $K = K(\varepsilon_2)$ we assume $\lambda=b\in K$ and obtain
the formula (1). The case 2) is clear because in this case $f(x)$ is
irreducible over $K$. In the case 3) we have again $\lambda = b\in K$ and
$e_0 = f_1 \in K^tG$, $e_{2^{s-1}} = f_2\in K^tG$. The other
idempotents in this case are obtained from (5) summing the pairs of
conjugates over$K$. In this way the formula (2) is completely
established. In the case 4)
by Lemma 2.4 we have $\lambda = b\varepsilon_{s+1}$, $b\in K$. The minimal
idempotents in this case are also obtained from (5)
summing the pairs of
conjugates over $K$ and this gives the formula (3).
\medskip

As a consequence of Theorem 3.1 one can obtain the minimal idempotents
of the factor-algebra $K[x]/I$, where $I$ is the ideal of the algebra
$K[x]$ generated by the polynomial $x^{2^n} - a$, $a\in K^{\ast}$. They are
obtained from the idempotents from Theorem 3.1 assuming that
$\bar g = x + I$.

  We shall note that when $\langle g\rangle$ is a cyclic 2-group and $K$
is a field of the first kind with respect to $2$ then the problem of finding
the explicit form of the minimal idempotents of the twisted group algebra
$K^t\langle g\rangle$ is open because of the serious problems that arrise
with the decompositions of the polynomials $x^{2^n}-a$ into irreducible
factors over the field $K$.

\vskip 1.truecm

\centerline{\bf REFERENCES}
\medskip

\item{1.} {A.\thinspace A. Bovdi,} Crossed products of semigroups and
rings,
{\sl Sib. Mat. Zh.} {\bf 4}, 1963, 481-500 (Russian).
\par

\item{2.} {Gr. Karpilovsky,} Field Theory, Marcel Dekker, Inc.,
New York and Basel, 1988.
\par

\item{3.} {N.\thinspace A. Nachev, T.\thinspace Zh. Mollov,} Minimal idempotents
of semisimple twisted group algebras of cyclic $p$-group of odd
order (Russian), {\sl Publ. Math., Debrecen} {\bf 35}, 1988, 309-319.
\par

\item{4.} {D.\thinspace S. Passman,} On the semisimplicity
of twisted group algebras, {\sl Proc. Amer. Math. Soc.} {\bf 25}, 1970,
161-166.
\par

\item{5.} {D.\thinspace S. Passman,}
The Algebraic Structure of Group Rings,
A Wiley-Interscience Publ., John Wiley and Sons, Inc., 1972.
\par

\item{6.} {B.\thinspace L. Van der Waerden,} Algebra I,
Sechste Auflage der Modernen
Algebra, Springer-Verlag, Berlin-G\"ottingen-Heidelberg, 1964.
\par

\bye